\documentclass[a4paper,12pt]{article}

\usepackage{hyperref}
\usepackage{stmaryrd}
\usepackage{amsfonts}
\usepackage{bbm}
\usepackage{amscd}
\usepackage{mathrsfs}
\usepackage{latexsym,amssymb,amsmath,amscd,amscd,amsthm,amsxtra,xypic}
\usepackage[dvips]{graphicx}
\usepackage[utf8]{inputenc}
\usepackage[T1]{fontenc}
\usepackage{lmodern}
\usepackage{nicefrac,mathtools,enumitem}
\usepackage{microtype}

\newtheorem{thm}{Theorem}[section]
\newtheorem{lem}[thm]{Lemma}
\newtheorem{cor}[thm]{Corollary}
\newtheorem{pro}[thm]{Proposition}

\newtheorem{rmk}[thm]{Remark}
\newtheorem{defi}[thm]{Definition}

\newcommand{\be }{\begin{equation}}
\newcommand{\ee }{\end{equation}}

\newcommand{\pf}{\begin{proof}}
\newcommand\Qed{\end{proof}}

\newcommand{\huaE}{\mathcal{E}}

\newcommand{\huaD}{\mathcal{D}}

\newcommand{\huaJ}{\mathcal{J}}

\newcommand{\frkg}{\mathfrak g}

\newcommand{\frkp}{\mathfrak p}

\newcommand{\frkX}{\mathfrak X}

\newcommand{\half}{\frac{1}{2}}

\newcommand{\pair}[1]{\left\langle #1\right\rangle}

\newcommand{\ppair}[1]{\left ( #1\right )_+}

\newcommand{\Courant}[1]{\left\llbracket  #1\right\rrbracket }

\newcommand{\br}[1]{   [ \cdot,    \cdot  ]   }

\begin{document}
\title{Twisted Courant algebroids and coisotropic Cartan geometries}

\author{Xiaomeng Xu \\
 Department of Mathematics and LMAM, Peking University, \\Beijing
100871, China\\\vspace{3mm} email: xuxiaomeng@pku.edu.cn}
\date{}
\maketitle{}

\begin{abstract} In this paper, we show that associated to any coisotropic Cartan geometry there is a twisted Courant algebroid. This includes in particular parabolic geometries. Using this twisted Courant structure, we give some new results about the Cartan curvature and the Weyl structure of a parabolic geometry. As more direct applications, we have Lie 2-algebra and 3D AKSZ sigma model with background associated to any coisotropic Cartan geometry.
\end{abstract}
\footnotetext{{\it{MSC}}: 53D17, 53A99.}
\footnotetext{{\it{Keyword}: Courant algebroids, Cartan geometry, parabolic geometry, Weyl structure, Pontryagin class}}

\section{Introduction}
Courant algebroids were introduced in \cite{lwx} to describe
the double of a Lie bialgebroid. They were further investigated in recent years, see
\cite{firstP}, \cite{DLB} for examples. Recently, the discovery in
\cite{4form} shows that the generalization of $4-$form twisted Courant algebroids arises naturally in the Courant sigma model with a Wess--Zumino term. The nature of this twist can also be explained by means of canonical functions and twisted QP manifolds, see \cite{IkedaXu}. At the same time, the further generalization of pre-Courant algebroids arise naturally in the construction of Courant algebroids (see \cite{pre}) and were further investigated in \cite{Pon}.  It unifies kinds of twisted Courant algebroids and the associated higher structures (see \cite{meixi}, \cite{Leibniz2}).

On the other hand, Cartan geometry is a common generalization of Riemannian geometry, conformal geometry and Klein geometry, which generalizes the linear tangent spaces of the former to the more general homogeneous spaces of the latter. It has been a unifying framework and a powerful tool for the study of conformal, projective, CR and related contact structures in differential geometry. Roughly speaking, a Cartan geometry is a curved analogon of a homogeneous space twisted by a curvature $\kappa$ (see \cite{3}, \cite{weylstructure}). Coisotropic Cartan geometry is a special case of Cartan geometry which includes parabolic geometry.
Pre-Courant algebroids were introduced in \cite{parageom} in this context. At the same time, any transitive quadratic Lie algebroid has a standard cotangent extension to a twisted Courant algebroid(\cite{dissection}), as a special example, one has twisted Courant structures from parabolic geometries(\cite{ams}).

In this paper, we deepen the connection between these two mathematical fields by using a more intrinsic way. Further, we can formulate some properties of the Cartan curvature and the Weyl structure of a parabolic geometry, because the associated twisted Courant structure includes much information of the original parabolic geometry. At the same time, Lie 2-algebras arise naturally as the underlying algebra structures of (pre-)Courant algebroids(\cite{Pon},\cite{roy2}), so we can get Lie 2-algebra structure from any coisotropic Cartan geometry where the graded space are cotangent bundle and tractor bundle. Also, one has topological field theories come from twisted Courant algebroids(\cite{4form}), so it is direct to obtain a 3D AKSZ sigma model with background associated to any coisotropic Cartan geometry. This sigma model is a generalization of 3D Chern-Simons theory, even if in the case of homogeneous space, it is also interesting and remains to be understood. We hope that further investigation of this twisted Courant structure and sigma model will give more insight into parabolic geometry or coisotropic geometry.

The paper is organized as follows. In section 2, we recall the definition and construction of pre-Courant algebroids, the Pontryagin class of a pre-Courant algebroid, and the criterions for a pre-Courant algebroid to be a twisted Courant algebroid.  In section 3, we review some basic knowledge about Cartan
geometry. Then we focus on coisotropic Cartan geometry and show the twisted Courant algebroid arising canonically from it. Using the properties of twisted Courant structures, we derive some properties of the Cartan curvature and Weyl structure of a parabolic geometry. In section 4, we give the Lie 2-algebra and the 3D AKSZ sigma model with background associated to any coisotropic Cartan geometry.

\section{Pre-Courant algebroids}
\begin{defi}\label{defi:preCA}
A pre-Courant algebroid consists of a vector bundle
$E\rightarrow{M}$, a fiberwise nondegenerate pseudo-metric
$\pair{\cdot,\cdot}$(by help of which we can identify $E$ with $E^*$), a bundle
map $\rho:E\rightarrow{TM}$ called the anchor and an $\mathbb
R-$bilinear bracket operation $\circ$ on $\Gamma(E)$ called Dorfman
bracket, such that for all $f\in{C^{\infty}(M)}$ and
$e_1,e_2,e_3\in{\Gamma(E)}$, the  following axioms are satisfied:
\begin{itemize}
\item[\rm(i)] $\rho({e_1}\circ{e_2})=[\rho(e_1),\rho(e_2)]$;
\item[\rm(ii)] $e_{1}\circ{e_{1}}=\half\partial\pair{e_{1},e_{1}}$;
\item[\rm(iii)]
$\rho(e_{1})\pair{e_{2},e_{3}}=\pair{e_{1}\circ{e_{2}},e_{3}}+\pair{e_{2},e_{1}\circ{e_{3}}}$,
\end{itemize}
where $\partial:\mathcal{C}^\infty(M)\rightarrow\Gamma(E)$ is the
$\mathbb R-$linear map defined by $$\pair{\partial f,e}=\rho(e)f.$$
\end{defi}

Following \cite{Pon}, for a pre-Courant algebroid $E$, the Jacobitor $J:\Gamma(E)\otimes\Gamma(E)\otimes\Gamma(E)\longrightarrow\Gamma(E)$ is defined by
\[J(e_1,e_2,e_3)=e_1\circ(e_2\circ e_3)-(e_1\circ e_2)\circ e_3-e_2\circ(e_1\circ e_3), \  \forall e_1,e_2,e_3\in\Gamma(E),\]
and the {\bf Pontryagin class} $\mathcal{P}:\Gamma(E)\times\Gamma(E)\times\Gamma(E)\times\Gamma(E)\mapsto{C^\infty}(M)$ is defined by
\[\mathcal{P}(e_1,e_2,e_3,e_4)=\pair{J(e_1,e_2,e_3),e_4}, \  \forall e_1,e_2,e_3,e_4\in\Gamma(E).\]

\begin{thm}[\cite{Pon}]\label{form}\rm
For any pre-Courant algebroid $E$ , we have $\mathcal{P}\in\Gamma(\wedge^4E)$ and $\mathcal{D}\mathcal{P}=0$, where $\mathcal{D}\mathcal{P}$ is defined by:
\begin{gather}
  \begin{split}
{\mathcal{D}\mathcal{P}(e_1,e_2,e_3,e_4,e_5)}=&\,\sum_{i=1}^4(-1)^{i+1}\rho(e_i)
(\mathcal{P}(e_1,\cdot\cdot\cdot,\hat{e_i},\cdot\cdot\cdot,e_5))  \\
 &+\sum_{i<j}(-1)^{i+j}\mathcal{P}((e_i\circ e_j),
e_1,\cdot\cdot\cdot,\hat{e_i},\cdot\cdot\cdot,\hat{e_j},\cdot\cdot\cdot,e_5).
  \end{split}\label{formulapartial}
\end{gather}
\end{thm}

$E$ is said to be a {\bf transitive pre-Courant algebroid}, if $\rho(E)=TM$. In this case, $(\ker\rho)^\perp$ is a subbundle of $E$, and $\Gamma((\ker\rho)^\perp)$ is an ideal of $\Gamma(E)$ under the Dorfman bracket, so we have an induced operation on $\Gamma(E/(\ker\rho)^\perp)$.
A {\bf transitive twisted Courant algebroid} is a transitive pre-Courant algebroid $E$ with a closed $4$-form $H$
such that $\mathcal{P}=\rho^*H$, where $\mathcal{P}$ is the Pontryagin class of $E$, i.e., for every $e_1,e_2,e_3\in\Gamma(E)$, we have
\begin{equation}\label{twistedCourant}
\rho^*(H(\rho(e_1),\rho(e_2),\rho(e_3),\rho(e_4))=\pair{e_1\circ (e_2\circ e_3)-(e_1\circ e_2)\circ e_3-e_2\circ(e_1\circ e_3),e_4}.
\end{equation}
If $H$ is zero, we recover the notion of Courant algebroid. We will
denote a transitive twisted Courant algebroid by
$(E,\pair{\cdot,\cdot},\circ,\rho,H)$.
As an application of Theorem \ref{form}, we have the following criteria for a pre-Courant algebroid to be a twisted Courant algebroid.

\begin{pro}\label{pro:strong}
Let  $(E,\pair{\cdot,\cdot},\circ,\rho)$ be a transitive pre-Courant algebroid and $\mathcal{P}$ be its Pontryagin class. Then the following statements are equivalent:
\begin{itemize}
\item[\rm(S1)]
  $J(e_1,e_2,e_3)\in{\Gamma((\ker\rho)^\perp)}$ for every
  ${e_1,e_2,e_3}\in{\Gamma(E)}$;

\item[\rm(S2)] $J(e,\cdot,\cdot)=0, $ for every $ e\in\Gamma(\ker\rho)$;

\item[\rm(S3)] $(E,\pair{\cdot,\cdot},\circ,\rho,H)$ is a twisted Courant algebroid, where $H$ is uniquely determined by $\mathcal{P}=\rho^*(H)$;
\item[\rm(S4)] $E/(\ker \rho)^\perp$ is a Lie algebroid with the induced operation is skew-symmetric.
\end{itemize}
\end{pro}
In order to construct a twisted Courant algebroid, we just need to start from a pre-Courant algebroid and to test one of the conditions listed in \ref{pro:strong}. On the other hand, without the algebra information given by $\circ$, a pre-Courant algebroid is just a pseudo-Euclidean vector bundle $(E,\pair{\cdot,\cdot})$ with a
bundle map $\rho: E\longrightarrow{TM}$, such that the image of
$\rho^*: T^*M\longrightarrow{E}$ ($E$ identified to $E^*$) is isotropic.
Such a triple $(E,\pair{\cdot,\cdot},\rho)$ is called
{\bf Courant vector bundle} by Izu Vaisman \cite{pre}. Summarizing the results in
\cite{pre}, we have the following method to construct pre-Courant algebroids:
\begin{pro}{\rm\cite{pre}}\label{pro:construct}
For any Courant vector bundle $E$, a pair $(\nabla,\beta)$ gives a
pre-Courant algebroid structures on $E$, where
$\nabla:\frkX(M)\times\Gamma(E)\longrightarrow\Gamma(E)$ is a metric
connection on $E$, i.e. $$X\pair{e_1,e_2}=\pair{\nabla
_Xe_1,e_2}+\pair{e_1,\nabla _Xe_2},$$
 and $\beta\in\Gamma(\wedge^2E^*\otimes{E})$
satisfying the following properties:
\begin{itemize}
\item[$(a)$] For any $e_1,e_2,e_3\in\Gamma(E)$, the map $(e_1,e_2,e_3)\longmapsto
\pair{\beta(e_1,e_2),e_3}$ is totally skew-symmetric with respect
to the pseudo-metric $\pair{\cdot,\cdot}$,

\item[$(b)$]
$\rho(\beta(e_1,e_2))=[\rho(e_1),\rho(e_2)]-\rho(\nabla_{\rho(e_1)}e_2-
\nabla_{\rho(e_2)}e_1).$
\end{itemize}
Using the pair $(\nabla,\beta)$, the pre-Courant algbroid bracket is
given by:
\[e_1\circ{e_2}=\nabla_{\rho(e_1)}e_2-
\nabla_{\rho(e_2)}e_1+\pair{\nabla
e_1,e_2}+\beta(e_1,e_2),\forall{e_1,e_2}\in\Gamma(E),\] where
$\pair{\nabla e_1,e_2}$ is defined by:
$$\pair{\pair{\nabla
e_1,e_2},e}=\pair{\nabla_{\rho(e)}e_1,e_2},\quad\forall{e}\in\Gamma(E).
$$
\end{pro}

\section{Twisted Courant algebroids from coisotropic Cartan geometry and applications }


We will construct a pre-Courant algebroid structure from coisotropic Cartan geometries.  This generalizes the case of parabolic geometry examined by Armstrong \cite{parageom}. Further, We will prove that these pre-Courant algebroids are twisted-Courant algebroids with the help of Proposition~\ref{pro:strong}.  Finally, we will give some applications in parabolic geometry.

At first, let us recall some basic definitions in Cartan geometry.
\begin{defi}\label{d:Cart}
Let $G$ be a Lie group with
Lie algebra $\frkg$, $P$ be a subgroup of $G$ with Lie algebra $\frkp$. A Cartan geometry $(\pi:\mathcal{G}\longrightarrow M,\omega,G/P)$ of type $(G,P)$ is
given by a principal bundle $\mathcal{G}\longrightarrow{M}$ with structure group $P$ equipped with a
$\frkg-$valued $1$-form $\omega$ (the Cartan connection) satisfying the following conditions:
\begin{itemize}
\item[\rm(a)]
the map $\omega_p:T_p\mathcal{G}\longrightarrow{\frkg}$ is a linear isomorphism for every
$p\in{\mathcal{G}}$;

\item[\rm(b)] ${R_a}^*\omega=Ad(a^{-1})\circ\omega$ $\forall{a}\in{P}$, where $R_a$ denotes the right action of the element $a$ in the structure group $P$;

\item[\rm(c)] $\omega(\zeta_{A})=A,$ $\forall{A}\in{\frkp},$ $\frkp$ is the Lie algebra of $P$ and $\zeta_{A}$ is the fundamental vector field on $\mathcal{G}$ generated by $A\in\frkg$.
\end{itemize}

\end{defi}

The {\bf curvature} of a Cartan connection is defined as
\begin{eqnarray}\label{curvature}
\kappa:=d\omega+\frac{1}{2}[\omega,\omega].
\end{eqnarray}
It is a $\frkg-$valued $2-$form $\kappa\in\Gamma(\wedge^2T\mathcal{G}\otimes\frkg)$.

For any homogeneous space $(G,P)$, with the left invariant Cartan $1$-form $\omega$, we see that they satisfy the conditions of Definition \ref{d:Cart} and the curvature $\kappa=d\omega+\frac{1}{2}[\omega,\omega]=0$ (Maurer--Cartan equation).  So homogeneous spaces provide flat Cartan geometries, and for a Cartan geometry $(\pi:\mathcal{G}\longrightarrow M,\omega,G/P)$, we call the corresponding homogeneous space $(G,P)$ its flat model.

One simple property we will use below is that the curvature is a
horizontal form, i.e. for $\xi\in\Gamma(\ker\pi_*)$, we have $\kappa(\xi,\cdot)=0$. Moreover, $\kappa$ is equivariant because of the
equivariance of $\omega$. Hence, we can view $\kappa$ as a $2-$form on $M$ with values in the vector bundle
$\mathcal{G}\times_P\frkg$ (the tractor bundle, see \cite{3}) associated to the adjoint representation of $P$ on $G$. This is to say, we can
view $\kappa$ as a function valued in $\wedge^2(\frkg/\frkp)^*\otimes \frkg$, i.e.
$\kappa\in{C^\infty(\mathcal{G},\wedge^2(\frkg/\frkp)^*\otimes \frkg)}$.

Starting from here, we consider parabolic geometry, i.e., Cartan geometry where $\frkg$ is semisimple and $\frkp$ is a parabolic
subalgebra of $\frkg$. We know that (see \cite{Pontryagin class}) for a semi-simple Lie algebra $\frkg$, a parabolic subalgebra $\frkp$ is
equivalent to having an $|l|-$grading of $\frkg$:
\[\frkg=\frkg_{-l}\oplus\cdot\cdot\cdot\oplus\frkg_{-1}\oplus\frkg_0\oplus\frkg_{1}\oplus\cdot\cdot\cdot\oplus\frkg_l,\]
where $\frkp=\frkg_0\oplus\frkg_{1}\oplus\cdot\cdot\cdot\oplus\frkg_l$, $[\frkg_i,\frkg_j]\subset\frkg_{i+j}$ and $\frkg_-=\frkg_{-l}\oplus\cdot\cdot\cdot\oplus\frkg_{-1}$ is generated by $\frkg_{-1}$$(\frkg_-$ can be identified with $\frkg/\frkp)$ and we have the following facts about this grading:

 $(a)$ There exists a grading element $E\in\frkg_0$ such that $[E,A]=iA$ if and only if $A\in\frkg_i$.

 $(b)$ The Killing form $B(\cdot,\cdot)$ defines isomorphisms $\frkg_{-i}\cong{\frkg^*_i}$ and $B(\frkg_i,\frkg_j)=0$ for all $j\neq-i$.

Now, the curvature of a parabolic geometry could be defined to be the curvature
function:

$\kappa:\mathcal{G}\longrightarrow \wedge^2\frkg_-^*\otimes \frkg$ \[p\longmapsto \kappa_p \ with \
\kappa_p(A,B)=d\omega_p(\omega_p^{-1}(A), \omega_p^{-1}(B))+[A,B], \ A,B\in \frkg_.\]

Let $(\pi:\mathcal{G}\longrightarrow M,\omega,G/P)$ be a Cartan geometry, we call it a {\bf coisotropic Cartan geometry}(following \cite{coisotropic}), if $\frkg$
is a quadratic Lie algebra and $\frkp$ is a coisotropic subalgebra of $\frkg$. In the case of a complex semisimple Lie algebra $\frkg$, $\frkp$ coisotropic implies that $\frkp$ is parabolic in $\frkg$ (see \cite{coisotropic}).  So in this case a cosiotropic Cartan geometry is just a parabolic geometry.  But in general, a quadratic Lie algebra may not be a semisimple Lie algebra and there exists a coisotropic Cartan geometry which is not a parabolic geometry.

The following theorem provides twisted-Courant algebroids from these rich geometric models.
\begin{thm}\label{para}
Let $(\pi:\mathcal{G}\longrightarrow M,\omega,G/P)$ be a coisotropic Cartan geometry, then there is canonically a twisted Courant algebroid structure on the tractor bundle $E=\mathcal{G}{\times_P}\frkg$.
\end{thm}
To prove this theorem, we need the following well-known result from Cartan geometry.
\begin{lem}\label{Liealgebroid}\cite{3}
Let $(\pi:\mathcal{G}\longrightarrow M,\omega,G/P)$ be a Cartan geometry, then there is a canonical Lie algebroid structure on $\mathcal{G}{\times_P}\frkg$, and this Lie algebroid is isomorphic to the Atiyah algebroid $T\mathcal{G}/P$.
\end{lem}
\pf  It is well-known that there is a bijection between sections of $\mathcal{G}\times_P\frkg$ and $P-$equivariant functions in $C^\infty(\mathcal{G},\frkg)$. For any $e\in\Gamma(\mathcal{G}\times_P\frkg)$, we denote $\widehat{e}$ the $P-$equivariant function corresponding to $e$. By the condition $(b)$ in the Definition 3.1, it is easy to see that $\omega^{-1}(\hat{e})$ is a $P-$invariant vector field on $\mathcal{G}$, and $P-$invariant vector fields on $\mathcal{G}$ form a Lie subalgebroid under the commutator bracket of vector fields (the Atiyah algebroid $T\mathcal{G}/P$), so there is a natural Lie bracket on $\Gamma(\mathcal{G}\times_P\frkg)$ which is defined by:
\[<e_1,e_2>:=\omega([\omega^{-1}(\hat{e_1}),\omega^{-1}(\hat{e_2})])=\nabla_{\rho(e_1)}e_2-
\nabla_{\rho(e_2)}e_1-[e_1,e_2]_\frkg-\kappa(e_1,e_2),\]
where $e_1,e_2\in\Gamma(E)$, $\nabla$ is the connection on $E$ induced by the Cartan connection(see \cite{3} for more details), $[\cdot,\cdot]$ is the Lie bracket of vector fields on $P$ and $[\cdot,\cdot]_\frkg$ is Lie bracket on $\Gamma(\mathcal{G}\times_P\frkg)$ induced by the pointwise Lie bracket on $\Gamma(\mathcal{G}\times\frkg)$.

On the other hand, we define a bundle map $\rho$ from $\mathcal{G}\times_P\frkg$ to $TM$ by
\begin{eqnarray*}
\rho(e_1)=\pi_*(\omega^{-1}(\hat{e_1})), \quad\forall \ e_1\in\Gamma(\mathcal{G}\times_P\frkg).
\end{eqnarray*}
It is easy to see that $\rho(<e_1,e_2>)=[\rho(e_1),\rho(e_2)]$ for every $e_1,e_2\in\Gamma(\mathcal{G}\times_P\frkg)$, so $(\mathcal{G}\times_P\frkg,<\cdot,\cdot>,\rho)$ is a Lie algebroid over $M$. 
\Qed

\pf[The proof of Theorem \ref{para}] The bundle $E=\mathcal{G}{\times_P}\frkg$ inherits a bilinear form $B(\cdot,\cdot)$ from the $P-$invariant
bilinear form on $\mathcal{G}\times{\frkg}$ and a Lie algebra bracket $[\cdot,\cdot]_\frkg$. So it is not hard to
see that $E$ is a Courant vector bundle with an anchor $\rho:E\longrightarrow{TM}$ given by
$\rho=\pi_*\circ\omega^{-1}$, where $\omega$ is the Cartan connection and $\pi$ is the natural projection
from $\mathcal{G}$ to $M$. Further, we have $\ker\rho=\mathcal{G}\times_P\frkp$ and $(\ker\rho)^\perp=\mathcal{G}\times_P\frkp^\perp$.

The next step is to choose a pair $(\nabla,\beta)$ satisfying the conditions in Proposition \ref{pro:construct}. Let $\nabla$ be the connection
on $E$ induced by the Cartan connection, we know $\nabla$ is a metric connection(see \cite{3} for more details). In order to
choose $\beta$, recall that there is a natural Lie algebroid structure on $\Gamma(E)$ whose Lie bracket is given by:
\[ <e_1,e_2>=\nabla_{\rho(e_1)}e_2-\nabla_{\rho(e_2)}e_1-[e_1,e_2]_\frkg-\kappa(e_1,e_2),\]
for every $e_1,e_2\in\Gamma(E)$.
Thus we have the following relation:
\[ [\rho(e_1),\rho(e_2)]=\rho<e_1,e_2>=\rho(\nabla_{\rho(e_1)}e_2-
\nabla_{\rho(e_2)}e_1))+\rho(-[e_1,e_2]_\frkg-\kappa(e_1,e_2)).\]
For every $e_1,e_2\in\Gamma(E)$, we define $B(e_1,\kappa(e_2,\cdot))\in\Gamma(E)$ by
\[B(B(e_2,\kappa(e_1,\cdot)),e_3)=B(e_2,\kappa(e_1,e_3)).\] Because $\kappa(e,\cdot)=0,$ for $\forall e\in ker\rho$, we
have $B(e_1,\kappa(e_2,\cdot))\in\Gamma((\ker\rho)^\perp)$, then $\rho\circ B(e_1,\kappa(e_2,\cdot))=0$.

From the discussion above and the conditions $(a)$ and $(b)$ of Proposition \ref{pro:construct}, we see that $\beta$ could be
defined by the skew-symmetric part of $-[e_1,e_2]_\frkg-\kappa(e_1,e_2)$:
\[\beta=-[e_1,e_2]_\frkg-\kappa(e_1,e_2)+B(e_2,\kappa(e_1,\cdot))-B(e_1,\kappa(e_2,\cdot)),\]
which satisfies the
conditions of Proposition \ref{pro:construct}, and the pre-Dorfman bracket on $E$
is given by:
\begin{eqnarray}\label{arm}
e_1\circ{e_2}=\nabla_{\rho(e_1)}e_2-
\nabla_{\rho(e_2)}e_1-[e_1,e_2]_\frkg+B(\nabla{e_1},e_2)-\kappa(e_1,e_2)+B(e_2,\kappa(e_1,\cdot))-B(e_1,\kappa(e_2,\cdot)),
\end{eqnarray}
where $B(\nabla e_1,e_2)\in\Gamma(E)$ is defined by
$B(B(\nabla e_1,e_2),e_3)=B(\nabla_{\rho(e_3)} e_1,e_2),$
for every $e_1,e_2,e_3\in\Gamma(E)$.
With the use of Lie bracket $<\cdot,\cdot>$, the above formula can be writen as:
\[e_1\circ{e_2}=<e_1,e_2>+B(e_2,\kappa(e_1,\cdot))-B(e_1,\kappa(e_2,\cdot))+B(\nabla{e_1},e_2).\]
Now we prove that the pre-Courant algebroid we get above is a twisted Courant algebroid.
Recall that $B(e_2,\kappa(e_1,\cdot))$,
$B(e_1,\kappa(e_2,\cdot))$ and $B(\nabla{e_1},e_2)$ are sections of subbundle $(\ker\rho)^\perp=\rho^*T^*M$($E^*$ and $E$ identified by $B(\cdot,\cdot)$) and the spaces of sections of $Ker\rho$ and ${\ker\rho}^\perp$ are two-sided ideals of $\Gamma(E)$, w.r.t.\ the Lie bracket $<\cdot,\cdot>$ on $\Gamma(E)$. So the operator $\circ$ induces a Lie bracket on sections of the bundle
$E/(\ker\rho)^\perp$. According to Proposition \ref{pro:strong}, $E$ is therefore a twisted-Courant
algebroid. \Qed

\begin{rmk}
According to Theorem \ref{form}, the Pontryagin class $\mathcal{P}$ of $E$ is a
section of $\wedge^4E$. On the other hand, the Cartan connection gives a bijection between
$\Gamma(E)$ and $P-$invariant vector fields on $\mathcal{G}$, so $\mathcal{P}$ could be seen as a $P-$invariant vertical
$4-$form on $\mathcal{G}$, and the formula $\mathcal{D}\mathcal{P}=0$ in Theorem \ref{form} is equivalent to say that $\mathcal{P}$ is a
closed $4-$form on $\mathcal{G}$. We will see that it is actually the first Pontryagin class of the quadratic Lie algebroid $E/(\ker\rho)^\perp$.
\end{rmk}

\begin{rmk}
It is easy to see that when the first Pontryagin class of the coisotropic Cartan geometry is zero, we can find a 3-form to twisted the bracket in \ref{arm} to get a Courant algebroid. The 3-form is unique up to an element in $H^3(M)$. See \cite{lx} for more details.
\end{rmk}
As an application of Theorem \ref{para}, we consider parabolic geometry which is coisotropic Cartan geometry obviously. Following the computation in \cite{parageom}, the Jacobiator of the twisted Courant algebroid given by Theorem \ref{para} associated to a parabolic geometry $(\pi:\mathcal{G}\longrightarrow M,\omega,G/P)$ is:
\begin{gather}\label{computation}
  \begin{split}
  J(e_1,e_2,e_3)=&\,-[e_1,B(e_3,\kappa(e_2,\cdot))-B(e_2,\kappa(e_3,\cdot))]_\frkg+B(e_3,\kappa(e_2,\kappa(e_1,\cdot))) \\
&-B(e_2,\kappa(e_3,\kappa(e_1,\cdot)))+c.p.,
  \end{split}
\end{gather}
where $e_1,e_2,e_3\in\Gamma(E=\mathcal{G}\times_P\frkg)$.

On the other hand, according to Proposition \ref{pro:strong}, for a transitive twisted Courant algebroid, $J(e,\cdot,\cdot)=0,$ for every $e\in{\ker\rho}=\mathcal{G}\times_P \frkp$, so we get a new result about the Cartan curvature in the following
corollary:
\begin{cor}
 For every parabolic geometry the following formula holds
\begin{gather}
  \begin{split}
 &[e_1,B(e_3,\kappa(e_2,\cdot))-B(e_2,\kappa(e_3,\cdot))]_\frkg-[e_3,B(e_1,\kappa(e_2,\cdot))]_\frkg\\
 &+[e_2,B(e_1,\kappa(e_3,\cdot))]_\frkg+B(e_1,\kappa(e_2,\kappa(e_3,\cdot))-\kappa(e_3,\kappa(e_2,\cdot)))\\
 &=0
  \end{split}
\end{gather} for every $e_1\in\Gamma(\mathcal{G}\times_P \frkp), \ e_2,e_3\in\Gamma(\mathcal{G}\times_P \frkg)$.\\
Especially,
\begin{equation} \label{formula}
  [e_1,B(e_2,\kappa(e_3,\cdot))]_\frkg=[e_2,B(e_1,\kappa(e_3,\cdot))]_\frkg,
\end{equation}
for every $e_1,e_2\in\Gamma(\mathcal{G}\times_P \frkp), \ e_3\in\Gamma(\mathcal{G}\times_P \frkg)$.
\end{cor}
Formula \eqref{formula} shows an important compatibility relation between the Cartan curvature of a parabolic geometry and the Lie algebra pair $(\frkg,\frkp)$ of its flat model.

To get further insight into the connection between parabolic geometry and the twisted Courant algebroid associated to it, we consider the Weyl structure of a parabolic geometry. Let $(\pi:\mathcal{G}\longrightarrow M,\omega,G/P)$ be a parabolic geometry and $\frkg$ have the grading $\frkg=\frkg_{-l}+\cdot\cdot\cdot+\frkg_{-1}+\frkg_0+\frkg_{1}+\cdot\cdot\cdot+\frkg_l$ corresponding to the parabolic subalgebra $\frkp$.  We can form an underlying principal $G_0-$bundle $\pi:\mathcal{G}_0\longrightarrow M$ by forming the quotient
$\mathcal{G}_0:=\mathcal{G}/P_+$, where $P_+$ is the nilpotent normal subgroup of $P$ with Lie algebra
$\frkp^\perp$ and $G_0$ is the Lie subgroup of $P$ with Lie algebra $\frkg_0$. We denote the quotient projection from $\mathcal{G}$ to $\mathcal{G}_0$ by $\pi_0$.
\begin{defi}{\rm\cite{weylstructure}}
Let
$(\pi:\mathcal{G}\longrightarrow{M},\omega,G/P)$ be a parabolic geometry. A global Weyl structure for
$\mathcal{G}$ is a global smooth $G_0-$equivariant section
$\sigma:\mathcal{G}_{0}\longrightarrow\mathcal{G}$ of the quotient projection $\pi_0$.
\end{defi}
\begin{pro}[\cite{weylstructure}]
There exists a global Weyl structure for every parabolic geometry.
\end{pro} A Weyl structure is an important aspect of parabolic geometry, see \cite{weylstructure} for
more details. We give the connections between the Weyl structure of a parabolic geometry with the twisted Courant algebroid associated to it. Firstly, recall that for any twisted Courant algebroid, the quotient bundle $E/(\ker\rho)^\perp$ has a natural Lie algebroid structure, called ample Lie algebroid of $E$, see \cite{dissection} for the details.
\begin{pro}\label{Ponclass}
The Atiyah algebroid $T\mathcal{G}_0/G_0$ is isomorphic to
Lie algebroid $E/(\ker\rho)^\perp$, where $E=\mathcal{G}{\times_P}\frkg$ is the twisted Courant algebroid given by Theorem \ref{para}.
\end{pro}
\pf Recall that any $e\in\Gamma(E)$ corresponds to a $P-$invariant vector field
$\widetilde{e}$ on $\mathcal{G}$, it is easy to see the quotient projection $\pi_0$ induces a bundle map
$\Phi$ from $E$ to $T\mathcal{G}_0/G_0$ by $e\longrightarrow{\pi_0}_*(\widetilde{e})$ , and the kernel of this
map is just $\mathcal{G}\times_P \frkp_+=(\ker\rho)^\perp$. So the bundle map induces a bundle isomorphism $\Phi'$ from $E/{(\ker\rho)}^\perp$ to $T\mathcal{G}_0/G_0$. On the other hand, when we consider the Lie algebra structure from the proof of Proposition \ref{Liealgebroid} and Theorem \ref{para}, we note that the bundle map $\Phi'$ is actually a Lie algebroid isomorphism. \Qed
\begin{rmk}
Suppose $E$ is the twisted Courant algebroid associated to a coisotropic Cartan geometry, thus $E/(\ker\rho)^{\perp}$ is just the Atiyah algbroid of the principal bundle $\pi:\mathcal{G}/P^+\longrightarrow M$, where $P_+$ is the nilpotent normal subgroup of $P$ with Lie algebra
$\frkp^\perp$. Following \cite{lx}, the Courant structure of $E$ could be realized by a reduction from the exact Courant algebroid $T\mathcal{G}_0\oplus T^*\mathcal{G}_0$. Further, it could also be realized by a reduction from the {\bf twisted action} algebroid (see \cite{Pon}) on $T\mathcal{G}\oplus T^*\mathcal{G}$ associated with the Cartan geometry $(\pi:\mathcal{G}\longrightarrow M,\omega,G/P)$. See \cite{lx} for more details about the exact realizations of (pre-)Courant algebroids.
\end{rmk}
Recall that for a twisted Courant algebroid $(E,\pair{\cdot,\cdot},\circ,\rho,H)$ on $M$, the first Pontryagin class $H\in\Omega^4(M)$ of the quadratic Lie algebroid $E/(\ker\rho)^\perp$ is given by
\[\rho^*H(e_1,e_2,e_3,e_4)=\pair{J(e_1,e_2,e_3),e_4}.\]
So equation \eqref{computation} and Proposition \ref{Ponclass} give a formula for
the first Pontryagin class of the principal $G_0-$bundle $\mathcal{G}_0$ on $M$ in terms of the Cartan curvature.
\begin{cor}\label{equality}
Let $(\pi:\mathcal{G}\longrightarrow M,\omega,G/P)$
be a parabolic geometry with Cartan curvature $\kappa$. Then the first Pontryagin class
$H\in\Omega^4(M)$ of the principal $G_0-$bundle $\mathcal{G}_0$ on $M$ is given by
\begin{gather}
  \begin{split}
  \pi^*H(e_1,e_2,e_3,e_4)=\,&-B([e_1,B(e_3,\kappa(e_2,e_4))-B(e_2,\kappa(e_3,e_1))]_\frkg,e_4)\\
  &+B(e_3,\kappa(e_2,\kappa(e_1,e_4)))-B(e_2,\kappa(e_3,\kappa(e_1,e_4))) \\
  &+c.p.,
  \end{split}
\end{gather}
where $e_1,e_2,e_3,e_4$ are right invariant vector fields on $T\mathcal{G}$.
\end{cor}
\begin{pro}
If $B(\kappa(e_1,e_2),e_3)$ is skew-symmetric for every $e_1,e_2,e_3\in\Gamma{E}$, then the first Pontryagin class $H\in\Omega^4(M)$ of the principal $G_0-$bundle $\mathcal{G}_0$ vanishes.
\end{pro}

\pf It is easy to see, with the same strategy as in the Theorem \ref{para}, the following equality defines a
pre-Courant structure $(E_1,\rho,g,\circ)$:
\[e_1\circ{e_2}=\nabla_{\rho(e_1)}e_2-\nabla_{\rho(e_2)}e_1+B(\nabla{e_1},e_2)-[e_1,e_2]_\frkg-\kappa(e_1,e_2).\] Using the same
computation as in \rm\cite{parageom}, we see
\[e_1\circ(e_2\circ e_3)-(e_1\circ{e_2})\circ e_3-e_2\circ(e_1\circ{e_2})=0,\]
that is to say $(E_1,\rho,B(\cdot,\cdot),\circ)$ is a Courant algebroid. By Proposition \rm\ref{Ponclass} the
quadratic Lie algebroid associated to $E_1$ is isomorphic to the Atiyah algebroid $T\mathcal{G}_0/G_0$.
On the other hand, the Pontryagin class of the ample Lie algebroid of a Courant algebroid is always zero (see \cite{firstP},\cite{dissection}).  So we obtain that the first
Pontryagin class of the principal $G_0-$bundle $\mathcal{G}_0$ is zero.\Qed

Corollary \ref{equality} gives the first Pontryagin class of the important principal bundle associated to a parabolic geometry in terms of the Cartan curvature.  In the case of a homogeneous space, the first Pontryagin class of the corresponding principal bundle vanishes.

\section{Lie 2-algebras associated to coisotropic Cartan geometries}
Following \cite{roy2} and \cite{Pon} , Lie 2-algebra arises naturally as the underlying algebra structure of a (pre-)Courant algebroid $E$, it encodes the Dorfman bracket and anchor operator on the $Z_2-$graded space given by cotangent bundle of base manifold and the Courant bundle $E$. Now, following the last section, for any coisotropic Cartan geometry $(\pi:\mathcal{G}\longrightarrow M,\omega,G/P)$ ,we have a twisted Courant algebroid on the tractor bundle $E=\mathcal{G}{\times_P}\frkg$, so we can get Lie 2-algebra structure on the graded space $(\Omega^1(M),E)$ naturally.

Firstly, for any (pre-)Courant algebroid $(E,\ppair{\cdot,\cdot},\rho,\circ)$, let us recall the Lie 2-algebra construction procedure below. Define the skew-symmetric bracket:
\begin{equation}
  \Courant{e_1,e_2}=\half(e_1\circ e_2-e_2\circ e_1)=e_1\circ
  e_2-\half\huaD\ppair{e_1,e_2}.
\end{equation}
Denote by $\huaJ:\wedge^3\Gamma(E)\longrightarrow \Gamma(E)$ its
Jacobiator:
$$
\huaJ(e_1,e_2,e_3)= \Courant{e_1,\Courant{e_2,e_3}}+c.p.
$$
By means of  \cite[Proposition 2.6.5]{roy2}, we have
$$
\huaJ(e_1,e_2,e_3)=J(e_1,e_2,e_3)-\huaD T(e_1,e_2,e_3),
$$
where $T(e_1,e_2,e_3)$ is given by
$$
T(e_1,e_2,e_3)=\frac{1}{6}\big(\ppair{\Courant{e_1,e_2},e_3}+c.p.\big).
$$

Now we define a 2-term complex as follows:
\begin{equation}\label{eqn:defihuaE}
  \huaE:\Omega^1(M)\stackrel{\rho^*}{\longrightarrow}\Gamma(E).
\end{equation}
Define degree-0 operation $l_2:\wedge^2\huaE\longrightarrow\huaE$ by
\begin{equation}\label{eqn:defil2}
  \left\{\begin{array}{rcll} l_2(e_1,e_2)&=& \Courant{e_1,e_2} & \mbox{in~ degree-0},~\forall~e_1,e_2\in\Gamma(E)\\
  l_2(e_1,\kappa)&=& \Courant{e_1,\kappa} &
  \mbox{in~
  degree-1},~\forall~e_1\in\Gamma(E),\kappa\in\Omega^1(M).
   \end{array}\right.
\end{equation}
Define degree-1 operator $l_3:\wedge^3\huaE\longrightarrow \huaE$ by
\begin{equation}\label{eqn:defil3}
 \begin{array}{rcll} l_3(e_1,e_2,e_3)&=& \huaJ(e_1,e_2,e_3), & \mbox{in~ degree-0},~\forall~e_1,e_2,e_3\in\Gamma(E).\\
   \end{array}
\end{equation}

\begin{pro}\cite{Pon}
  For a twisted Courant algebroid $E$, $(\huaE,l_2,l_3)$ is a Lie
  2-algebra, where $\huaE,l_2,l_3$ are given by
  \eqref{eqn:defihuaE},
  \eqref{eqn:defil2}, \eqref{eqn:defil3}
  respectively.
\end{pro}

So eventually for any coistropic Cartan geomtry $(\pi:\mathcal{G}\longrightarrow M,\omega,G/P)$, we get a 2-term complex as follows:
\begin{equation}
  \huaE:\Omega^1(M)\stackrel{\rho^*}{\longrightarrow}\Gamma(E=\mathcal{G}{\times_P}\frkg).
\end{equation}
Define a degree-0 operation $l_2:\wedge^2\huaE\longrightarrow\huaE$ by
\begin{eqnarray*}
l_2(e_1,e_2)&=&\nabla_{\rho(e_1)}e_2-
\nabla_{\rho(e_2)}e_1-[e_1,e_2]_\frkg-\kappa(e_1,e_2)+B(e_2,\kappa(e_1,\cdot))-B(e_1,\kappa(e_2,\cdot))\\
&&+\frac{1}{2}B(\nabla{e_1},e_2)-\frac{1}{2}B(\nabla{e_2},e_1).
\end{eqnarray*}
Define degree-1 operator $l_3:\wedge^3\huaE\longrightarrow \huaE$ by
\begin{equation}\label{eqn:defil3}
 \begin{array}{rcll} &&l_3(e_1,e_2,e_3)\\
 &=& -[e_1,B(e_3,\kappa(e_2,\cdot))-B(e_2,\kappa(e_3,\cdot))]_\frkg+B(e_3,\kappa(e_2,\kappa(e_1,\cdot)))
-B(e_2,\kappa(e_3,\kappa(e_1,\cdot)))\\&&-\huaD\{\frac{1}{6}(B(\nabla_{\rho(e_1)}e_2,e_3)-B(\nabla_{\rho(e_2)}e_1,e_3)+\frac{1}{2}B(\nabla_{\rho(e_3)}e_1,e_2)-\frac{1}{2}B(\nabla_{\rho(e_3)}e_2,e_1))
\\&&+ \frac{1}{6}(B([e_1,e_2]_\frkg,e_3)-B(\kappa(e_1,e_2),e_3)+B(e_2,\kappa(e_1,e_3))-B(e_1,\kappa(e_2,e_3)))\}+c.p..
   \end{array}
\end{equation}
\begin{cor}
$(\huaE,l_2,l_3)$ is a Lie
  2-algebra, where $\huaE,l_2,l_3$ are given by the equations above.
\end{cor}

\section{Further Outlook}
In this very first step, we combine the Courant algebroid with the coisotropic Cartan geometry. We see that coisotropic Cartan geometries are abundance sources of twisted Courant algebroids. One the other hand, the tool developed in Courant algebroid is used to get new results about coisotropic Cartan geometry.

In order to obtain more insight into parabolic geometry or coisotropic Cartan geometry from the view of mathematical physics, a reasonable way is to combine topological field theory with it. This is based on the fact that 3D AKSZ sigma models with background one-to-one correspond to twisted Courant algebroids, see \cite{4form} and \cite{royl}. By bulk-boundary corresponding, this is equivalent to a 4D AKSZ sigma model with boundary. Now in the case of target space with a coisotropic Cartan structure, we apply exactly the same methods to the twisted Courant algebroids from coisotropic geometry, then we get a topological field theory associated to any coisotropic Cartan structure where we have a 4-from background on target given by the first Pontryagin class. In a local chart, the classical action functional is given explicitly by a formula involving the Lie algebra, Killing form and Cartan curvature. Even in the case of homogeneous space, this sigma model generalizing Chern-Simons theory remains to be understood. It is naive to hope that the 3D AKSZ sigma model will open a new door to coisotropic Cartan geometries.

\end{document}